\documentclass[12pt,a4paper,reqno]{amsart}
\usepackage{amsmath}
\usepackage{amsfonts}
\usepackage{amssymb}
\numberwithin{equation}{section}

\theoremstyle{plain}

\newtheorem{theorem}[subsection]{Theorem}

\newtheorem{lemma}[subsection]{Lemma}

\theoremstyle{definition}

\newtheorem{definition}[subsection]{Definition}

\renewcommand{\leq}{\leqslant}
\renewcommand{\geq}{\geqslant}

\newsavebox{\proofbox}
\savebox{\proofbox}{\begin{picture}(7,7)%
  \put(0,0){\framebox(7,7){}}\end{picture}}
\def\boxeq{\tag*{\usebox{\proofbox}}}

\def\endproof{\hfill{\usebox{\proofbox}}}

\def\Z{\mathbb{Z}}

\def\ni{\noindent}
\def\vs{\vspace{11pt}}

\def\gen{{\rm gen\,}}
\def\diam{{\rm diam\,}}

\begin{document}

\title{Sets with small sumset and rectification}

\author{Ben Green}
\address{Department of Mathematics\\
University of Bristol\\
University Walk\\
Bristol BS8 1TW\\
England
}
\email{b.j.green@bristol.ac.uk}

\author{Imre Z. Ruzsa}
\address{Alfr\'ed R\'enyi Institute\\ 
Budapest\\
}
\email{ruzsa@renyi.hu}

\thanks{While this work was carried out the first author was resident in Budapest, and was supported by the \textit{Mathematics in Information Society} project carried out by the R\'enyi Institute, in the framework of the European Community's \textit{Confirming the International R\^ole of Community Research} programme. The second author was supported by the Hungarian National Foundation for Scientific Research
(OTKA), Grants No. T 29759 and T 38396.}

\begin{abstract}
We study the extent to which sets $A \subseteq \mathbb{Z}/N\mathbb{Z}$, $N$ prime, resemble sets of integers from the additive point of view (``up to Freiman isomorphism''). We give a direct proof of a result of Freiman, namely that if $|A + A| \leq K|A|$ and $|A| < c(K)N$ then $A$ is Freiman isomorphic to a set of integers. Because we avoid appealing to Freiman's structure theorem, we get a reasonable bound: we can take $c(K) \geq (32K)^{-12K^2}$.\vs

\ni As a byproduct of our argument we obtain a sharpening of the second author's result on sets with small sumset in torsion groups. For example if $A \subseteq \mathbb{F}_2^n$, and if $|A + A| \leq K|A|$, then $A$ is contained in a coset of a subspace of size no more than $K^22^{2K^2 - 2}|A|$.
\end{abstract}

\maketitle

\section{Introduction} \label{sec1}
\ni Additive questions involving integers are often simpler than the analogous
questions for residues modulo a prime. For instance, the conjecture of Erd\H os
and Heilbronn \cite{e-h} that given $k$ residues modulo a prime $N$, the number of sums of
distinct pairs is at least $ \min (2k-3,N)$ was open for over 30 years until
Dias da Silva and Hamidoune proved it in \cite{dds-h}, while the corresponding statement for
integers is obvious.\vs

\ni By \emph{rectification} we mean a result that allows us to replace a set
of residues by a set of integers which behaves identically from our point of
view. To illuminate this terminology we quote Freiman \cite{freiman73}: ``It
is known that sets of integers and sets of residues are different if
considered from an additive viewpoint: addition of sets of integers takes
place on ``the line'' whereas addition of sets of residues is performed
``on the circle''. However in the case of noniterated addition of sets of
residues which are not too large, this circle is ``broken up''. This
phenomenon expresses itself by the fact that such a set of residues is
isomorphic to the corresponding set of integers.''\vs

\ni To make this more exact, let $A,B$ be subsets of abelian groups, and let $f : A \rightarrow B$ be 1-1.
We say that $f$ is a \emph{Freiman isomorphism of order
$k$}, or $F_k$-isomorphism for short, if for $a_1, \dots , a_k, a_1', \dots , a_k'
\in A$ the equation \[ a_1 + \dots  + a_k = a_1' + \dots  + a_k'\]
 holds if and only if
\[ f(a_1) + \dots  + f(a_k) = f(a_1') + \dots  + f(a_k').\] This means that additive
properties involving at most $k$ summands are preserved under such a map. For
the important case $k=2$ we shall use the term Freiman-isomorphism. We say that
a set of residues is \emph{rectifiable of order $k$}, if it is $F_k$-isomorphic
to a set of integers (simply \emph{rectifiable}, if $k=2$).\vs

\ni Suppose our set of residues modulo a prime $N$ is contained in
$\{m, m+1, \dots , m+l\}$. The map which takes this interval into the interval
$\{m, \dots , m+l\}$ of integers is a $F_k$-isomorphism if $kl<N$, and so such a set
is rectifiable of order $k$ if $l<N/k$. Since multiplication by a nonzero
residue is clearly a Freiman isomorphism of any order, the same conclusion
holds for a set contained in an arithmetic progression of length $<N/k$. Using
the box principle one can show that this happens whenever the set has
$< \log N/ \log k $ elements, and we cannot say more if we only know the
number of elements.\vs

\ni The situation changes if we also know that our set has few sums. We write
     \[   A + B = \{a+b: a\in A, b\in B \},  \]
     $A-B$ is defined similarly, and we denote the $k$-fold iterated sum by
     \[  kA = A + \dots  + A  .   \]
     Now Freiman \cite[Ch. III, \S 4]{freiman73} proved that for any $K$ there is a constant $c(K) > 0$ with the following property: if $A \subseteq \Z/N\Z$, if $|A + A| \leq K|A|$ and if $|A| < c(K)N$,
then $A$ is rectifiable. A simpler proof can be found in \cite{r98c}.
Both proofs apply Freiman's
fundamental description of sets of integers with small \textit{doubling constant}, that is to say satisfying $|A + A| \leq K|A|$ where $K$ is to be thought of as very small in comparison with $|A|$. One formulation of this result asserts that
such a set is contained in a generalized arithmetic progression
     \[   P(a; q_1, \dots , q_d; l_1, \dots , l_d) = \left \{a + \sum_{i=1}^d  x_i q_i: 0 \leq 
x_i\leq l_i-1   \right \},  \]
     where the \emph{dimension} $d$ of $P$ is bounded by $d\leq f_1(K)$ and its
\emph{size} $|P| = l_1 \dots  l_d$ (which may be equal to its cardinality or may
be larger) satisfies $|P| \leq  f_2(K)n$. For different versions and proofs of this
result see Freiman \cite{freiman73}, Ruzsa \cite{r92i,r93c}, Bilu \cite{bilu99}
and Chang \cite{chang03}; the best estimates of $f_1, f_2$ are given in Chang's
work.\vs

\ni The aim of this paper is to give a direct proof of the rectification
theorem. 

     \begin{definition}
     The \emph{diameter} $\diam A$ of a set $A$ \textup{(}in  ${\mathbb {Z}}$ or ${\mathbb {Z}}/m\mathbb{Z}$\textup{)} is defined as the
smallest integer $l$ for which there exists some $a,d$ such that $A\subset \{a, a+d, \dots , a+ld\}$.
     \end{definition}

\ni Here are our main results.  In these statements, $\alpha \in (0,1)$ and $K \geq 1$ are arbitrary, whilst $k \geq 2$ is an integer.                      

     \begin{theorem}\label{thm1}
     Let $N$ be a prime and let $A\subseteq {\mathbb {Z}}/N\mathbb{Z}$ be a set with $|A|=\alpha N$ and $\min (|2A|,|A - A|) = K|A|$. Suppose that $\alpha \leq (16 K)^{-12K^2}$.
Then the diameter of $A$ is at most 
     \begin{equation} \label{diambound}
     12\alpha^{1/4K^2}\sqrt{\log(1/\alpha)}N.
     \end{equation}
          \end{theorem}

     \begin{theorem}\label{thm2}
     Let $N$ be a prime and let $A\subseteq {\mathbb {Z}}/N\mathbb{Z}$ be a set with $|A|=\alpha N$ and $|2A|=K|A|$. Then $A$
is rectifiable of order $k$ provided that $\alpha \leq (16 k K)^{-12K^2}$.
     
     \end{theorem}

     \begin{theorem}\label{thm3}
     Let $A\subseteq {\mathbb {Z}}$ be a set with $|A|=n$ and $|2A|=Kn$. Then $A$
is $F_k$-isomorphic to a set contained in the interval $[1, (16 k K)^{12K^2}n ]$.
     
     \end{theorem}

\ni As a byproduct, we will obtain an improved version of the analog of
Freiman's theorem in torsion group (see Section 6).\vs

\ni In \S \ref{sec2}, we deduce a covering property for sets with small doubling
and use it to estimate iterated sums. In \S \ref{sec3} we connect diameter
and large Fourier coefficients. In \S \ref{sec4} we deduce the existence of large
Fourier coefficients from the covering property. \S \ref{sec5} contains the
synthesis.\vs

\section{Doubling, covering and the size of iterated sums}\label{sec2}

\ni The \emph{covering property} of a set means that its sumset is covered by
a few translations of the original set. We will show that if $A$ has a small
sumset, then its difference set possesses such a property, and this will 
be used to find bounds for the cardinality of its higher order sumsets.

     \begin{lemma}
     Let $A, B_1, B_2$ be finite sets in an arbitrary commutative group. Write
$|A|=n$, $|A+B_i|=K_i n$. There is a set $T\subseteq B_1+B_2$ such that $|T|\leq 2K_1 K_2-1$
and
     \begin{equation} \label{inc}
     B_1-B_1+B_2-B_2 \subseteq  A-A+T-T .  
     \end{equation}
     \end{lemma}

\begin{proof}
\ni For the proof we will use the following Pl\"unnecke-type inequality of the
second author which may be found in \cite{r89e}, Theorem 6.1:

     \begin{lemma}
     Let $A, B_1, B_2$ be finite sets in an arbitrary commutative group. Write
$|A|=n$, $|A+B_i|=K _i n$. There is a set $A' \subseteq  A $ such that 
     \[   |A' +B_1+B_2 | \leq  K _1 K_2 |A'| .   \]
     \end{lemma}

\ni Now take this set $A'$ and select elements $t_1, t_2, \dots , t_k$ of
$B_1+B_2$, as many as possible, with the property that
     \[  \left |  (A'+t_j) \setminus  \bigcup _{i=1}^{j-1} (A'+t_i) \right | \geq  |A'|/2 . \]
     This property immediately implies that
     \[  \left |  \bigcup _{i=1}^{k} (A'+t_i) \right | \geq  (1+k/2) |A'| , \]
     and since this set is contained in $A'+B_1+B_2$, we can conclude that $k
\leq  2K_1K_2 -1$.\vs

\ni We set $T= \{t_1, \dots , t_k \}$ and we show that it has property \eqref 
{inc}. Take a generic element of $   B_1-B_1+B_2-B_2$, say $u_1 - u_2$ with
$u_j \in B_1 + B_2$. Since neither $u_j$ can serve as $t_{k+1}$, we have
     \[  \left |  (A'+u_j) \setminus  \bigcup _{i=1}^{k} (A'+t_i) \right | < |A'|/2 . \]
     This means that there is an $a\in A'$ such that
     \[    a+u_j \in  \bigcup _{i=1}^{k} (A'+t_i)  \]
     for $j=1,2$, say $a+u_j = a_j + t_{i_j}$. By substracting these equations
we obtain
     \[   u_1-u_2 = a_1-a_2 + t_{i_1} - t_{i_2} \in  A'-A'+T-T \subseteq  A-A+T-T \]
     as required.
     \end{proof}

     \begin{lemma} \label{Incm}
     Let $A$ be a finite set in an arbitrary commutative group. Write  $|A|=n$,
and assume that $ \min (|2A|, |A-A|) \leq  Kn$.
There is a set $T\subseteq 2A$ such that $|T|\leq 2K^2-1$
and such that for every positive integer $m$ we have
     \begin{equation} \label{incm}
     (m+1)(A-A) \subseteq  A-A+m(T-T) .
     \end{equation}
     \end{lemma}

     \begin{proof}
     We obtain the case $m=1$ by substituting $B_1=B_2=A$ or $B_1=B_2=-A$ into
\eqref{inc}. The general case follows by an obvious induction.
     \end{proof}\vs

\ni We think that it is often better to work with such a covering property
than with the small doubling assumption. (One reason is that the covering
property is preserved under a linear mapping, whereas the small doubling may not be.)

     \begin{definition}
     We say that a set $B$ \textup{(}in a group\textup{)} is \emph{$k$-covering} if there is a
set $T$ such that $|T|=k$ and  $B+B\subseteq B+(T-T)$ \textup{(}whence $(m+1)B\subseteq B+m(T-T)$ for every
positive integer $m$\textup{)}.
     \end{definition}

\ni This particular form of the inclusion comes from our above results. A more
complicated inclusion, say $B+B\subseteq B+100T-200T$ would be almost as useful, since
in estimating the size of $mB$ for large $m$ it is more important to have a
bound on the size of $T$ than on the number of its occurrences.\vs

     \ni We define $J(k,m) $ as the number of $k$-tuples of integers $(x_1, \dots ,
x_k)$ satisfying 
     \begin{equation} \label{a}
      \sum_{i=1}^k  x_i^+ = \sum_{i=1}^k  x_i^- \leq  m .
     \end{equation}

     \begin{lemma} \label{Estjcov}
     Let $B$ be a finite set in an arbitrary commutative group. Assume that $B$
is $k$-covering. For every positive integer $m$ we have
     \begin{equation} \label{estjcov}
     |(m+1)B| \leq  |B| J(k, m).
     \end{equation}
     \end{lemma}

     \ni\textit{Proof.} 
     Indeed, the inclusion $(m+1)B\subseteq B+m(T-T)$ yields
     \[  |(m+1) B| \leq  |B| |m(T-T)| .  \]
     If $T=\{t_1, \dots , t_k\}$, then every element of $m(T-T)$ can be written in
the form $\sum  x_i t_i$ with $x_i$ satisfying \eqref  a, hence
     \begin{equation}\boxeq  |m(T-T)| \leq  J(k, m). \end{equation}

     \begin{lemma}
     For $m\geq k$ we have
     \begin{equation} \label{b}
     J(k, m) < (14m/k)^k .
     \end{equation}
     \end{lemma}

     \ni\textit{Remark.} For large $m/k$ $J(k,m)$ behaves essentially as $(4em/k)^k$, so at most
the constant 14 can be improved.\vs

     \noindent\textit{Proof. }
     Write $x_i=y_i s_i$ with $y_i=|x_i|$ and $s_i = \pm 1$. The number of
possible choices for $s_i$ is $2^k$.  The $y_i$'s satisfy $\sum  y_i \leq  2m$, and the
number of such $k$-tuples is
     \[   \binom{2m+k}{k} < \frac{ (2m+k)^{2m+k}}{k^k (2m)^{2m}} \leq  \left ( 7m/k \right )^k. \]
     To see the last inequality here observe that with the notation $k=tm$,
$0\leq t\leq 1$ it reduces to
     \[   (2+t)^{2+t} \leq  4 \cdot 7^t .  \]
     This holds for $t=0$ and $t=1$, so it holds in between by the convexity of
the function $x \log x$. Thus we obtain
     \begin{equation}\boxeq   J(k,m) \leq  2^k \binom{2m+k}{k} \leq  (14m/k)^k. \end{equation}

     \ni By combining the previous two lemmas we obtain the following one.

     \begin{lemma} \label{Estecov}
     Let $B$ be a finite set in an arbitrary commutative group. Assume that $B$
is $k$-covering. For every positive integer $m\geq k$ we have
     \begin{equation} \label{estecov}
     |(m+1)B| < (14m/k)^{k} |B| .
     \end{equation}
     \end{lemma}

     \section{Large Fourier coefficients and diameter}\label{sec3}

    \ni Our plan is to use the growth rate estimate of the previous section to
find a large Fourier coefficient and to use it to estimate the diameter. This
cannot be achieved directly, since the union of an interval and a few random numbers
has a very large Fourier coefficient and a large diameter as well. We circumvent this
difficulty as follows. We show that a large Fourier coefficient implies that a high
percentage of our set is contained in a short progression, apply this to the
difference set, and show that if a large part of $A-A$ is in a progression,
then the whole set $A$ can be covered by a progression of the same length. This section
contains the tools needed to carry through this project.

     \begin{lemma}  \label{Cover}
     Let $A\subseteq  {\mathbb {Z}}/N\mathbb{Z}$ be a set of size $n$. Assume that
     $$  | (A - A)\setminus [b, b+l ]| < n/2   $$
with some $b$ and some $l<N/3$. Then $A\subseteq [a, a+l]$ for some $a$.
     \end{lemma}

     \begin{proof}
     As every quantity involved is invariant under shifts of $A$, we may assume
that $0\in A$ (whence $A\subseteq A - A$) and that the longest gap in $A$ is of the form $[k, -1]$.
We have $  | A\setminus [b, b+l ]| < n/2   $, so $A' = A \cap  [b, b+l ]$ satisfies
$|A'|>n/2$. Now take an arbitrary $a\in A$. The set $A'-a\subseteq A - A$ cannot be contained
in $(A - A)\setminus [b, b+l ]$, so $a'-a\in [b, b+l ]$ for some $a'\in A'$. Hence
     \[   a \in  a'- [b, b+l ] \subseteq  [-l, l ] .\]
     This holds for every $a$, so we have proved that $A\subseteq [-l , l]$.\vs

  \ni   This inclusion implies that there is a gap of length $(N- l) - l  >  l $.
As we assumed that the longest gap ends at 0, there can be no element of $A$ in
$[- l ,0)$ and the above inclusion can be strengthened to $A\subseteq [0, l  ]$.
     \end{proof}\vs

  \ni  Now let $B \subseteq \mathbb{Z}/N\mathbb{Z}$. For $r \in \mathbb{Z}/N\mathbb{Z}$ we define the Fourier transform of $B$ at $r$ by 
     \[   \widehat{B}(r) = \sum _{b\in B} e(br/N) \]
     where, as usual, $e(\theta) = e^{2\pi i \theta}$. Clearly $\widehat{B}(0) = |B|$ and $\sup_{r \neq 0 } |\widehat{B}(r)| \leq |B|$. We will be interested in situations where equality almost holds in the latter bound.

     \begin{lemma}
     Let $B\subseteq {\mathbb {Z}}/N\mathbb{Z}$. Assume that $\varepsilon ,\delta $ are real
numbers such that $0<\varepsilon <1$, $0<\delta <1/2$ and $|\widehat{B}(1)|\geq (1-2\varepsilon (1- \cos \pi \delta ))|B|$. Then
     \[   |B\setminus [a, a+l]| < \varepsilon |B|  \]
     with suitable $a$ and $l<\delta N$.
     \end{lemma}

  \ni   This is Lev \cite[Theorem 1]{lev04} with some changes in the notation. A
more comfortable sufficient assumption for the conclusion is
     \begin{equation} \label{lev}
     |\widehat{B}(1)|\geq (1-8\varepsilon \delta ^2)|B| ;
     \end{equation}
     this follows from Lev's sharp bound via the inequality
     \[   \cos \delta \pi  < 1-4\delta ^2 .  \]


     \begin{lemma} \label{Diam}
     Let $A\subseteq  {\mathbb {Z}}/N\mathbb{Z}$. Write $D=A-A$, and suppose that
     \[   |\widehat{D}(r)| \geq |D| - 4\delta^2|A| \]
     for some  $\delta \in (0, 1/3)$ and $r\ne 0$. Then
$\diam A < \delta N $.
     \end{lemma}

     \begin{proof} Write $|A| = n$ and $|D| = m$.
     Assume first that $r=1$. We apply the previous lemma (in the form \eqref  {lev})
     with $D$ in the place of
$B$ and $\varepsilon =n/2m$. We get an interval of length $<\delta N$ that contains $D$ with
less than $\varepsilon m=n/2$ exceptions. Now Lemma \ref{Cover} implies that $A$ is
contained in an interval of length $<\delta N$.\vs

  \ni    The general case follows by applying the special case to a suitable
dilation of $A$.
     \end{proof}

     \section{Finding a large Fourier coefficient}\label{sec4}

 \ni    With a future application in mind \cite{freiman-groups} we will describe the results here
for general finite abelian groups. Let $G$ be a finite Abelian group and let $B\subseteq G$. The Fourier transform $\widehat{B}$ is defined for characters $\gamma \in G^{\ast}$ by
     \[  \widehat{B}(\gamma ) = \sum _{b\in B} \gamma (b) . \]
     We denote the principal character by $\gamma _0$.

     \begin{lemma} \label{Largecoeff}
     Let $G$ be a finite Abelian group with cardinality $N$, and suppose that $B\subseteq G$ has size $\beta N$. Assume that
$B$ is $k$-covering for some integer $k\geq 2$, and that $\beta \leq 14^{-k - 1}$. Then there is a
character $\gamma \ne \gamma _0$ such that $|\widehat{B}(\gamma )| \geq  (1-\eta )|B|$, where
     \begin{equation} \label{largecoeff}
     \eta  = 18\beta^{1/k}\log(1/\beta)/k.
     \end{equation}
     \end{lemma}

     \begin{proof} Let $m$ be a positive integer, and write $R = |(m+1)B|$. 
     We will calculate high moments of $\widehat{B}$. We have
     \[   \widehat{B}(\gamma )^{m+1} = \sum _{x\in G} r_{m+1}(x) \gamma (x), \]
     where $r_{m+1}(x)$ the number of $(m+1)$-tuples $(b_1,\dots,b_{m+1}) \in B^{m+1}$ such that $b_1 + \dots + b_{m+1} = x$ (equivalently, $r_{m+1}$ is the $(m+1)$-fold autoconvolution of the characteristic function of $B$). We have
     \[   \sum_{x \in G}  r_{m+1}(x) =  \widehat{B}(\gamma _0)^{m+1} = |B|^{m+1} \]
     and by Cauchy-Schwarz we get
     \[   \sum_{x \in G}  r_{m+1}(x)^2 \geq   |B|^{2m+2}/R . \]
     Now Parseval's identity yields
     \[   \sum _{\gamma \in G^{\ast}}  |\widehat{B}(\gamma )|^{2m+2} = N \sum_{x \in G}  r_{m+1}(x)^2 \geq  N|B|^{2m+2}/R . \]
     On substracting the contribution of $\widehat{B}(\gamma _0)=|B|$ we obtain
     \[   \sum _{\gamma \ne \gamma _0} |\widehat{B}(\gamma )|^{2m+2} \geq  (N/R - 1)|B|^{2m+2}. \]
     Since
     \[   \sum _\gamma  |\widehat{B}(\gamma )|^{2} = N|B|, \]
     we can conclude that
     \[   \max _{\gamma \ne \gamma _0} |\widehat{B}(\gamma )|^{2m} \geq  (1/R - 1/N)|B|^{2m+1}.
\]

   \ni  A similar argument is used by Schoen \cite{schoen03}.\vs

   \ni  We use the estimate $R\leq  (14m/k)^k |B|$ from \eqref  {estecov} (valid
for $m\geq k$) and take ($2m$)-th roots to obtain
     \begin{equation} \label{c}
     \max _{\gamma \ne \gamma _0} |\widehat{B}(\gamma )| \geq  \left ( (k/148 m)^k - \beta \right )^{1/2m} |B| .
     \end{equation}
     This estimate is positive if $(k/14m)^k>\beta$, a condition which is satisfied by taking
     \[   m = \left\lfloor k(2\beta)^{-1/k}/14\right\rfloor  . \]
     The assumption $\beta \leq 14^{-k-1}$ implies that $m\geq k$. With this choice \eqref  c gives
     \begin{equation} \label{d}
     \max _{\gamma \ne \gamma _0} |\widehat{B}(\gamma )| \geq \beta^{1/2m} |B|.
     \end{equation}
     One easily sees that
     \[   m \geq  k\beta^{-1/k}/36  \]
     so that
     \[ \beta^{1/2m} = e^{-\log(1/\beta)/2m} > 1 - \frac{\log(1/\beta)}{2m} \geq 1 - \eta.\]
     Substituting this into \eqref  d gives \eqref  {largecoeff} immediately.
     \end{proof}

     \begin{lemma} \label{Largecoeff2}
     Let $G$ be a finite Abelian group of size $N$. Let $A \subseteq G$ be a set for which
     $ \min (|2A|, |A-A|) \leq  K|A|$. Write $D = A - A$, and suppose that $|D| = \tau N$ where $\tau \leq 14^{-2K^2}$. Then there is a character $\gamma \neq \gamma_0$ such that $|\widehat{D}(\gamma)| \geq (1 - \eta)|D|$, where
   \begin{equation} \label{largecoeff2}   \eta = 9K^{-2}\tau^{1/2K^2}\log(1/\tau).
     \end{equation}
     \end{lemma}

     \begin{proof}
     This follows from the previous lemma, using that $D$ is $k$-covering with
some $k\leq 2K^2-1$.
     \end{proof}

     \section{Assembling the pieces}\label{sec5}

  \ni   In this section we prove Theorems \ref{thm1}, \ref{thm2} and \ref{thm3}.\vspace{11pt}

    \noindent\textit{Proof of Theorem \ref{thm1}.}
      Write $D=A-A$. We have $|D| \leq K^2 |A|$ by a result of the second author
\cite{r78a}. In the notation of Lemma \ref{Largecoeff2}, this means that $\tau \leq K^2\alpha$. By Lemma
\ref{Largecoeff2} we know that there is an $r\ne 0$ such that $|\widehat{D}(r)|\geq (1-\eta )|D|$, the value of $\eta $ being given by \eqref  {largecoeff2} (recall that $\tau = |D|/N$).
By Lemma \ref{Diam} we know that 
$ \diam A < \delta N$ if $\eta |D| \leq 4\delta^2|A|$ and $\delta <1/3$. The first of these conditions will be satisfied if $\delta = 2K\eta^{1/2}$. The second will too, if we invoke the assumption $\alpha \leq (16 K)^{-12K^2}$ subject to which our theorem is being proved (this is a somewhat tedious check, substituting for $\eta $ from \eqref  {largecoeff2}). Another irritating calculation of an almost identical kind confirms that $\delta$ satisfies the bound claimed in Theorem \ref{thm1}.\endproof\vspace{11pt}

     \noindent\textit{Proof of Theorem \ref{thm2}.}
     We only have to check that the assupmtion on $\alpha$ implies that the $\delta$ of
Theorem \ref{thm1} satisfies $\delta <1/k$. This is straightforward.
     \endproof\vspace{11pt}

     \noindent\textit{Proof of Theorem \ref{thm3}.}
     With no loss of generality we may assume that $A$ is one of the shortest sets in its $F_k$
isomorphy class. In other words, $A\subseteq [1, L]$ and there is no isomorphic set contained
in $[1, L-1]$. Take a prime $N\in (kL, 2kL)$. Let $A'\subseteq {\mathbb {Z}}/N\mathbb{Z}$ be the image of $A$; since
$N > kL$, $A$ and $A'$ are $F_k$-isomorphic. Write $|A'| = \alpha N$. Let $l$ be the diameter of $A'$; we
already know that $l\leq L$. By definition, this means that there is an $A''\subseteq {\mathbb {Z}}/N\mathbb{Z}$ which is
homothetic (and a fortiori $F_k$-isomorphic) to $A'$ and contained in 
$\{1, 2,\dots , l\}$. The corresponding set $A'''$ of integers is a set isomorphic
to $A$ and contained in $[1, l]$; by the minimality of $L$ we infer $l=L$.
Comparing this with the bound of Theorem \ref{thm1} for $l$ we obtain
     \[ N \leq 2kL = 2kl \leq 24k\alpha^{1/4K^2}\sqrt{\log (1/\alpha)}N.\]
     This implies that $\alpha > (16kK)^{-12K^2}$, which means that $L \leq (16kK)^{12K^2}n$. This is what we wanted to prove.\endproof

     \section{Sets with small sumset in a torsion group}

\ni We recall that the \emph{exponent} of a group $G$ is
the smallest positive integer $r$ such that $rg=0$ for every $g\in G$.

     \begin{theorem}
     Let $G$ be a commutative group of exponent $r$, and assume that $A\subseteq G$.
     \item{(a)} {If $|2A|=K|A|$, then $A$ is contained in a coset of a subgroup
of size
     \[   \leq  K ^2 r^{\lfloor 2K ^2-2\rfloor}|A| . \] }
     \item{(b)} {If $|A-A|=K n$, then $A$ is contained in a coset of a subgroup
of size
     \[   \leq  K  r^{\lfloor 2K ^2-2\rfloor}|A| . \] }
     \end{theorem}

 \ni    The first result of this type is due to the second author \cite{r99b}.
In that paper the exponent of $r$ is $K ^4$. A different approach was developed by
Deshouillers, Hennecart and Plagne \cite{deshetal03}.
In its present forms, it works only for $r=2$ and $K <4$. For $K <3$ their
result is superior to ours.\vs

     \begin{proof} 
     For a set $X\subseteq G$, we write $\gen X$ to denote the
subgroup generated by $X$. Our set $A$ is contained in a coset of $\gen (A-A)$,
which is the union of the increasing sequence $m(A-A)$. Hence \eqref  {incm}
implies
          \[   \gen (A-A) \subseteq  (A-A) + \gen (T-T) . \]
          with a suitable set $T$, $|T|=k\leq 2K ^2-1$.
     If $T = \{t_1, \dots , t_k \}$, then $\gen (T-T)$ is generated by the
$k-1$ elements $t_i-t_1$, $i=2, \dots , k$. So if every element has order $\leq r$,
then
     \[   | \gen (T-T) | \leq  r^{k-1} \leq  r^{\lfloor 2K ^2-2\rfloor} ,\]
     consequently
     \[ | \gen (A-A) | \leq  | A-A| |\gen (T-T)| \leq  |A-A| r^{\lfloor 2K ^2-2\rfloor} . \]
     This proves part (b). To obtain part (a) we recall that if $|2A|=K n$, then
$|A-A|\leq K ^2n$ by a result of the second author \cite{r78a}.
     \end{proof}


\providecommand{\bysame}{\leavevmode\hbox to3em{\hrulefill}\thinspace}
\providecommand{\MR}{\relax\ifhmode\unskip\space\fi MR }
\providecommand{\MRhref}[2]{%
  \href{http://www.ams.org/mathscinet-getitem?mr=#1}{#2}
}
\providecommand{\href}[2]{#2}

     \end{document}